\documentclass[ejs]{imsart}

\usepackage{amsthm}
\usepackage{amsmath}
\usepackage{amssymb}
\usepackage{color}
\usepackage[draft]{graphicx} 
\graphicspath{{figures/}} 
\usepackage{wrapfig} 
\usepackage[font=normalsize,labelfont=bf]{caption} 
\usepackage{float} 
\usepackage{booktabs} 
\usepackage{tikz} 
\usepackage{hyperref} 
\usepackage[mathscr]{eucal} 

\doi{https://doi.org/10.1214/21-EJS1964}
\pubyear{2022}
\volume{16}
\firstpage{862}
\lastpage{891}

\newcommand{\R}{{\mathbb R}}
\newcommand{\N}{{\mathbb N}}

\newcommand{\Var}{{\rm Var}}

\newcommand{\Cov}{\operatorname{Cov}}
\newcommand{\F}{{\mathcal F}}
\newcommand{\claw}{\stackrel{\mathcal{D}}{\longrightarrow}}
\newcommand{\cpr}{\stackrel{\mathcal{P}}{\longrightarrow}}

\newtheorem{theorem}{Theorem}
\newtheorem{lemma}{Lemma}[section]
\newtheorem{definition}[lemma]{Definition}
\newtheorem{corollary}[lemma]{Corollary}
\newtheorem{proposition}[lemma]{Proposition}
\newtheorem{remark}{Remark}

\begin{document}

\begin{frontmatter}

\title{Change-point detection based on weighted two-sample U-statistics\support{H. Dehling and K. Vuk were supported by the Collaborative Research Grant SFB 823 \textit{Statistical modelling of nonlinear dynamic processes}. M. Wendler was supported by the German Research Foundation (DFG), project WE 5988/3 \textit{Analyse funktionaler Daten ohne Dimensionsreduktion.}
}}

\runtitle{Weighted change-point tests}


\begin{aug}
\author{\fnms{Herold} \snm{Dehling}\ead[label=e1]{herold.dehling@ruhr-uni-bochum.de}}
\and
\author{\fnms{Kata} \snm{Vuk}\ead[label=e2]{kata.vuk@ruhr-uni-bochum.de}}

\address{Fakult\"at f\"ur Mathematik, \\ Ruhr-Universit\"at Bo\-chum, \\ 44780 Bochum, Germany \\ \printead{e1,e2}}

\author{\fnms{Martin} \snm{Wendler}\ead[label=e3]{martin.wendler@ovgu.de}}
\address{Institut f\"ur Mathematische Stochastik, \\ Otto-von-Guericke-Universit\"at Magdeburg, \\ 39106 Magdeburg, Germany\\ \printead{e3}}

\runauthor{H. Dehling et al.}

\end{aug}

%
%
%
%

\begin{abstract}
We investigate the large-sample behavior of change-point tests based on weighted two-sample U-statistics, in the case of short-range dependent data. Under some mild mixing conditions, we establish convergence of the test statistic to an extreme value distribution. A simulation study shows that the weighted tests are superior to the non-weighted versions when the change-point occurs near the boundary of the time interval, while they loose power in the center.
\end{abstract}



\tableofcontents

\end{frontmatter}


\section{Introduction}
In this paper, we study nonparametric tests for change-points in time series that are based on weighted two-sample $U$-statistics. By a suitable choice of the weights, we obtain tests that are able to detect changes that occur very early or very late during the observation period. Our results cover both the CUSUM test and the Wilcoxon change-point test, as well as many other robust and non-robust tests. We investigate the large-sample behavior of our tests in the case of short-range dependent data under mild conditions that cover, e.g. ARMA and ARCH processes. By means of a simulation study, we analyze the small sample behavior of our tests, e.g. regarding robustness and the ability to detect early and late changes. As an application, we analyze a data set of daily stock returns of Wirecard during the weeks prior to the detection of accounting fraud in June 2020.

There is a vast body of literature on change-point tests using U-statistics. Gombay \cite{G:2000}, and Kirch and Stoehr \cite{KS:2020} investigate sequential change-point tests based on U-statistics. Wang, Volgushev and Shao \cite{WVS:2019} apply certain U-statistics for change-point detection in high-dimensional time series. Zhang, Wang, and Shao \cite{ZWS:2021} propose an adaptive change-point test based on U-statistics. Ra{\v{c}}kauskas and Wendler \cite{RW:2020} apply U-statistics for the detection of epidemic changes.

We assume that the data are generated by a stochastic process $(X_i)_{i\geq 1}$ which follows the model
\[
  X_i=\mu_i+\xi_i, \; i\geq 1, 
\]
where $(\mu_i)_{i\geq 1}$ is an unknown signal, and where $(\xi_i)_{i\geq 1}$ is a short-range dependent stationary stochastic process. Given the observations $X_1,\ldots,X_n$, we want to test the hypothesis that the signal is constant, i.e.
\[
  H_0\;:\; \mu_1=\ldots=\mu_n 
\]
against the alternative of a change in the mean at an unknown point $k^\ast$ in time, i.e.
\[
  H_1\; :\; \mu_1=\ldots =\mu_{k^\ast} \neq \mu_{k^\ast+1}=\ldots = \mu_n, \mbox{ for some } k^\ast \in \{1,\ldots,n-1\}.
\]
If the change-point $k^\ast$ was known, we would have a two-sample problem with the samples $X_1,\ldots,X_{k^\ast}$ and $X_{k^\ast+1},\ldots, X_n$, respectively, and we could apply standard tests such as the two-sample Student t-test or the Wilcoxon two-sample test for a change in mean. Up to normalization, both are special cases of two-sample $U$-statistics $\sum_{i=1}^{k^\ast} \sum_{j=k^\ast+1}^n h(X_i,X_j)$, with a suitably chosen kernel function $h:\R^2\rightarrow \R$. 

In the change-point setting, where a change occurs at an unknown point in time, we have a family of two-sample problems, indexed by the potential change-point $k$, and thus we are naturally led to the two-sample $U$-statistic process 
\[
  \sum_{i=1}^k \sum_{j=k+1}^n h(X_i,X_j), \; 0\leq k\leq n. 
\]
A variety of change-point tests can be derived from this process by taking suitable functionals such as weighted maxima 
\[
  \max_{1\leq k\leq n-1} \frac{1}{\big(\frac{k}{n}(1-\frac{k}{n})  \big)^\gamma } \frac{1}{n^{3/2}} 
  \Big| \sum_{i=1}^k \sum_{j=k+1}^n h(X_i,X_j)  \Big|,
\]
where $0\leq \gamma \leq \frac{1}{2}$ is some parameter to be chosen. 

For $\gamma=0$, we obtain non-weighted tests, which have been widely studied in the literature, starting with Darkhovskh \cite{D:1976} and Pettitt \cite{P:1979}, who studied the special case of a Wilcoxon-type test statistic. Cs\"org\H{o} and Horvath \cite{CsH:1988} investigated U-statistics with general kernels, in the case of independent data. They could show that the asymptotic distribution under the null hypothesis is the Kolmogorov-Smirnov distribution, which is the distribution of the supremum of a Brownian bridge.  Dehling, Fried, Garcia and Wendler \cite{DFGW:2015} extended these results to weakly dependent data.  Under long-range dependence, the limiting distribution is given by the supremum of a linear combination of Hermite processes; see Dehling, Rooch and Taqqu \cite{DRT:2013} for the Wilcoxon test, and Dehling, Rooch and Wendler \cite{DRW:2017} for arbitrary kernels. For $0<\gamma<\frac{1}{2}$, the limit distribution under independence is the supremum of the appropriately weighted Brownian bridge, see e.g. the seminal monograph by Cs\"org\H{o} and Horv\'ath \cite{CsH:1997}. The moment conditions for such a limit theorem have been relaxed by Cs\"org\H{o}, Szyszkowicz, and Wang \cite{CsSW:2008}.

In the present paper, we focus on the extreme case $\gamma=\frac{1}{2}$, where we obtain the test statistic 
\begin{equation}
T_n:=\max_{1\leq k\leq n-1} \frac{1}{ \sqrt{k(n-k)n  } } \Big| \sum_{i=1}^k \sum_{j=k+1}^n h(X_i,X_j)  \Big|
\end{equation}
Under the null hypothesis, after some suitable normalization and centering, $T_n$ converges in distribution to the Gumbel extreme value distribution. This has been derived by Cs\"org\H{o} and Horv\'ath \cite{CsH:1988} under independence. We will show that the same holds in the case of short-range dependent data. Antoch, Hu\v{s}kov\'{a} und Pr\'{a}\v{s}kov\'{a} \cite{AHP:1997} have studied the large-sample behavior of weighted versions of the CUSUM test for dependent observations, in particular for linear processes. We also show that the test is consistent against a wide class of alternatives. For independent data, the behavior under the alternative has been studied by Ferger \cite{F:1994} and Gombay \cite{G:2000}. We have conducted an extensive simulation study comparing this test with the corresponding non-weighted test. Our simulations confirm the intuition that the weighted tests have more power against very early and very late changes, while the non-weighted tests are more powerful against changes in the middle of the observation period. 

By the choice of the kernel function $h$, the weighted two-sample U-statistics lead to a flexible class of change-point tests. As special examples, we obtain the CUSUM test for $h(x,y)=y-x$, and the Wilcoxon test for $h(x,y)=1_{\{ x\leq y \}} -\frac{1}{2}$. More generally, one can take the kernels $h(x,y)=\psi(y-x)$ for some anti-symmetric function $\psi:\R \longrightarrow \R$. Depending on the choice of $\psi$, one obtains tests with specific properties, such as robustness against outliers, and tests that are powerful for certain alternatives.

The rest of the paper is organized as follows. In the next section, we present the detailed technical assumptions, and we give the main theoretical results of our paper. In section~3, we present the outcomes of a major simulation study comparing weighted and non-weighted tests as well as the robust Wilcoxon test and the non-robust CUSUM test. Full details of the proofs are presented in the final section.

\section{Main theoretical results}

In this section, we analyze the large sample behavior of the suitably normalized and centered test statistic $T_n$ under the hypothesis, and under a wide class of alternatives. A major ingredient in the proof is a
new  Darling-Erd\H{o}s type limit  theorem for the tied-down random walk of dependent random variables, 
which might also be of independent interest. Before we present our results, we give some definitions. Throughout this paper, the stochastic process $(X_i)_{i\geq 1}$ will be assumed to be $\alpha$-mixing in the sense of Rosenblatt \cite{R:1956}. 

\begin{definition}
The stochastic process $(X_i)_{i\geq 1}$ is said to be $\alpha$-mixing if
\[
  \alpha(k) := \sup_{n\geq 1} \sup \big\{ |P(A\cap B)-P(A)\, P(B)| : A\in \F_1^n, B\in \F_{n+k}^\infty\big\} \longrightarrow 0,
\]
as $k \longrightarrow \infty$, where $\F_a^b$ denotes the $\sigma$-field generated by the random variables $X_a,\ldots,X_b$. We define the generalized inverse $\alpha^{-1}:[0,1]\longrightarrow \N$ by 
\[
  \alpha^{-1}(u):= \min \big\{k\in \N: \alpha(k) \leq u \big\} = \sum_{i=0}^\infty 1_{ \{ \alpha(i) > u  \}  }.
\]
\end{definition}

Our theoretical results require assumptions on the rate of decay of the mixing coefficients $(\alpha(k))_{k\geq 1}$. We formulate these assumptions using a concept introduced by Rio \cite{Rio:1993} that is based on the quantile function which we define below.

\begin{definition}
For a random variable $X$, the upper quantile function $Q_X:[0,1]\longrightarrow \R$ is defined by 
\[
  Q_X(u)=\inf \big\{  t\in \R: P(X >t)\leq u   \big\}
\]

\end{definition}

Finally, we will assume that the kernel function $h:\R^2\rightarrow \R$ satisfies the variation condition, which is a continuity assumption introduced by Denker and Keller \cite{DK:1986}, and that the kernel has uniform $(2+\delta)$-moments.

\begin{definition}
A kernel $h:\R^2\rightarrow \R$ satisfies the variation condition if there exist constants $L>0$ and $\varepsilon_0>0$ such that for all $\varepsilon \in (0,\varepsilon_0)$ 
\begin{equation}
E \Big(  \sup_{\{(x,y): \|(x,y)-(X,Y)\| \leq \varepsilon \}}  |h(x,y)-h(X,Y)| \Big) \leq L\, \varepsilon,
\label{eq:var-cond}
\end{equation}
where $X,Y$ are independent random variables with the same distribution as $X_1$, and where $\|\cdot \|$ denotes the Euclidean norm on $\R^2$. 
\end{definition}

\begin{definition} 
Let $(X_i)_{i\geq 1}$ be a stationary process. A kernel $h$ has uniform $(2+\delta)$-moments if for all $k\in \mathbb{N}_0$ 
\begin{align*}
E|h(X_1,X_k)|^{(2+\delta)} &\leq M, \\
E|h(X,Y)|^{(2+\delta)} &\leq M,
\end{align*}
where $X,Y$ are independent copies of $X_1$, and where $M$ is a constant.
\end{definition}
The following theorem is the main theoretical result of this paper. Throughout, $Q_{|X|}$ will denote the common quantile function of the $X_k$'s.  

\begin{theorem}
\label{thm:mainresult}
Let $(X_i)_{i\geq 1}$ be an $\alpha$-mixing strictly stationary process. Let $h(x,y)$ be a bounded anti-symmetric kernel with uniform $(2+\delta)$-moments and satisfying the variation condition \eqref{eq:var-cond}. Moreover, assume that there exist constants $p>2$ and $\varepsilon>6/\delta$ such that
\begin{equation}
\int_0^1 (\alpha^{-1}(u))^{4+\varepsilon} Q_{|X|}^p(u) du <\infty.
\label{ass.epsilon}
\end{equation}
Then, under the null hypothesis $H_0$, as $n\rightarrow \infty$,
\[
  \frac{\sqrt{2\log\log n}}{\sigma_h} T_n -b_n \stackrel{\mathcal{D}}{\longrightarrow} G_2,
\]
where $G_2$ is the Gumbel extreme value distribution with distribution function 
\begin{equation}
 G_2(x) =\exp(-2\exp(-x)),
\label{eq:gumbel-2}
\end{equation}
and where the centering constants $b_n$ and the long-run variance $\sigma_h^2$ are defined by
\begin{align}
 b_n&=2\log\log n +\frac{1}{2} \log\log\log n -\frac{1}{2} \log \pi  \label{eq:b-n} \\
 \sigma_h^2&= \Var(h_1(X_1) +2 \sum_{k=2}^\infty \Cov(h_1(X_1),h_1(X_k)).
 \label{eq:lrv-h1}
\end{align}
Here, $h_1(x)$ denotes the first order term in the Hoeffding decomposition of $h$, as defined below. 
\end{theorem}

The idea of the proof is to apply the Hoeffding decomposition, which was introduced by Hoeffding \cite{H:1948}, and to show that the degenerate part is asymptotically negligible. Thus, it will remain to show that the linear part converges to the extreme value distribution $G_2$.  For a two-sample U-statistic, the Hoeffding decomposition of the kernel $h$ is given by
\begin{equation}
h(x,y) = \theta+h_1(x) +h_2(y) +\Psi(x,y), \label{hoeffding}
\end{equation}
where the terms on the right hand side are defined by
\begin{align*}
 \theta&=Eh(X,Y) \\
 h_1(x)&=Eh(x,Y)-\theta \\
 h_2(y)&=Eh(X,y)-\theta \\
 \Psi(x,y)&= h(x,y) -h_1(x)-h_2(y) -\theta,
\end{align*}
and where $X,Y$ are two independent random variables with the same distribution as $X_1$. Note that in our case, since $h$ is assumed to be anti-symmetric, $\theta=0$ and $h_2(x)=-h_1(x)$. Applying the Hoeffding decomposition of the kernel $h$ to the test statistic $T_n$, we obtain
\begin{align}
& \frac{1}{ \sqrt{ k(n-k)n } } \Big| \sum_{i=1}^k \sum_{j=k+1}^n h(X_i,X_j)  \Big| \label{decomp} \\
  = \, & \frac{1}{ \sqrt{ k(n-k)n } } \Big| \sum_{i=1}^k \sum_{j=k+1}^n \big(h_1(X_i)-h_1(X_j) +\Psi(X_i,X_j)\big)  \Big| 
 \nonumber \\
 =\, &  \frac{1}{ \sqrt{ k(n-k)n } } \Big| (n-k) \sum_{i=1}^k  h_1(X_i) -k \sum_{j=k+1}^n h_1(X_j) + \sum_{i=1}^k \sum_{j=k+1}^n \Psi(X_i,X_j)  \Big| 
 \nonumber \\
 =\, &  \Big|  \sqrt{\frac{n}{ k(n-k)}} \Big( \! \sum_{i=1}^k  h_1(X_i) -\frac{k}{n}  \sum_{j=1}^n h_1(X_j)\! \Big)  + 
  \frac{1}{ \sqrt{k(n-k)n} }\sum_{i=1}^k \sum_{j=k+1}^n \! \Psi(X_i,X_j)  \Big| 
 \nonumber 
 \end{align}
In order to show that $\frac{\sqrt{\log\log n}}{\sigma_h} T_n-b_n \stackrel{\mathcal{D}}{\longrightarrow} G_2$, it thus suffices to show that
\begin{equation}
 \sqrt{ \log\log n }  \max_{1\leq k\leq n-1} \frac{1}{\sqrt{ k(n-k)n }} \Big| \sum_{i=1}^k \sum_{j=k+1}^n \Psi(X_i,X_j)  \Big| \stackrel{P}{\longrightarrow} 0,
 \label{eq:deg-part}
\end{equation}
and that 
\begin{equation}
  \frac{\sqrt{2\log\log n}}{\sigma_h} \max_{1\leq k\leq n-1} \sqrt{ \frac{n}{k(n-k)} } \Big( \sum_{i=1}^k  h_1(X_i) -\frac{k}{n}  \sum_{j=1}^n h_1(X_j)\Big) \stackrel{\mathcal{D}}{\longrightarrow} G_2.
 \label{eq:lin-part}
\end{equation}
The asymptotic negligibility of the degenerate part, i.e. \eqref{eq:deg-part}, will be established in Proposition~5.1, while \eqref{eq:lin-part}  will be a consequence of a suitable Darling-Erd\H{o}s theorem; see Theorem~3 below.

In the next theorem, we investigate the large sample behavior of $T_n$ under the alternative $H_1$. We define
\[
  \Delta=\Delta_n:=Eh(X_1,X_n^\prime),
\]
where $X_n^\prime$ is independent of $X_1$ and has the same distribution as $X_n$.  Note that $\Delta_n$ measures the size of the change in the distribution of $X_i$ at the change point. 

\begin{theorem}
\label{thm:consistency}
Assume that the degenerate kernel $h$ has uniform $(2+\delta)$-moments and that it satisfies the variation condition \eqref{eq:var-cond}. Let 
\begin{equation}
\sum_{k=1}^\infty k \left( \int_{0}^{\alpha(k)} Q_{|X|}(u) du \right)^{\frac{\delta}{3+2\delta}} <\infty.
\label{eq:mix-2}
\end{equation}
Moreover, assume that the alternative $H_1$ holds, and that 
\begin{equation}
 \sqrt{ \frac{ k_n^\ast (n-k_n^\ast) }{ n \log \log n }  }  |\Delta_n| \longrightarrow \infty,
\label{eq:h1-local}
\end{equation}
where $k^\ast=k_n^\ast$ denotes the location of the change. Then
\begin{equation}
 \frac{1}{\sqrt{\log \log n}} T_n \stackrel{P}{\longrightarrow }\infty.
\end{equation}
\end{theorem}

\begin{corollary}
The test that rejects the null hypothesis $H_0$ of no change when 
\[
\frac{\sqrt{2\log\log n}}{\sigma_h} T_n -b_n  \geq g_{2,\alpha},
\]
where $g_{2,\alpha}$ denotes the upper $\alpha$-th quantile of the Gumbel distribution $G_2$, has asymptotic level $\alpha$. Moreover, the test is consistent against any alternatives $H_1$ that satisfy \eqref{eq:h1-local}.
\end{corollary}
\begin{proof}
Under the null hypothesis, the test statistic converges by Theorem~1 to the Gumbel extreme value distribution $G_2$, and thus the test that rejects the null hypothesis when the statistic exceeds $g_{2,\alpha}$ has asymptotically level $\alpha$. Regarding the behavior under the alternative, 
we will show that under the assumptions of the corollary
\begin{equation}
\frac{ \sqrt{2\log\log n}  }{ \sigma_h} 
T_n -b_n \stackrel{P}{\longrightarrow} \infty. 
\label{eq:tn-infty}
\end{equation}
Let $K>0$ be a given constant, then
\begin{align*}
 & P\Big( \frac{ \sqrt{2\log\log n}  }{ \sigma_h} 
T_n -b_n \geq K\Big) \\
= \, &  P \Big(  T_n  \geq \frac{(b_n+K)\sigma_h}{\sqrt{2\log\log n}} \Big) \\
= \, &  P \Big( \frac{1}{\sqrt{\log\log n}}  T_n  \geq \frac{(b_n+K)\sigma_h/\sqrt{2}}{\log\log n } \Big) \\
\geq \, & P \Big( \frac{1}{\sqrt{\log\log n}}  T_n  \geq 3 \sigma_h/\sqrt{2} \Big),
\end{align*}
for all $n$ large enough, since $(b_n+K)/\log\log n\rightarrow 2$ as $n\rightarrow \infty$. Now, the right hand side converges to $1$ by Theorem~2.
\end{proof}
\begin{remark}
The condition \eqref{eq:h1-local} puts restrictions on the time $k_n^\ast$ as well as the magnitude  $\Delta_n$ of the change, in order for our test to be consistent. For an early change, i.e. when $k_n^\ast = o(n)$, the condition
\eqref{eq:h1-local}  is equivalent to
\[
  \frac{k_n^\ast \Delta_n^2}{\log\log n}\rightarrow \infty,
\]
stating that $k_n^\ast \cdot \Delta_n^2$ has to grow faster than $\log\log n$ for the weighted test to be consistent. Keeping the magnitude of the change $\Delta_n\equiv \Delta$ constant,  this means that $k_n^\ast$ must grow faster than $\log\log n$. 
\end{remark}

In the next theorem, we derive the asymptotic distribution of the weighted tied-down random walk. This result will be used in the analysis of the asymptotic distribution of the linear part of the test statistic $T_n$. 

\begin{theorem}
\label{thm:tied-down}
Let $(X_i)_{i\geq 1}$ be an $\alpha$-mixing stationary process, let $S_n:=\sum_{i=1}^n X_i$ denote the partial sum process, and let 
\begin{equation}
\sigma^2=\Var(X_1) +2\sum_{k=2}^\infty \Cov(X_1,X_k)
\label{eq:lrv}
\end{equation}
denote the long-run variance.
If there exists a $2<p\leq 3$ such that
\begin{equation}
 \int_0^1 (\alpha^{-1}(u))^{p-1} Q_{|X|}^p (u) \, du<\infty,
\label{eq:mixing-3}
\end{equation}
then 
\[
 \frac{ \sqrt{2 \log\log n}}{\sigma} \max_{1\leq k\leq n-1} \sqrt{ \frac{n}{k(n-k)}}
  \Big| S_k-\frac{k}{n} S_n \Big| -b_n \stackrel{\mathcal{D}}{\longrightarrow}G_2,
\]
where $G_2(x)=\exp(-2\exp(-x))$, and where $b_n$ is defined as in \eqref{eq:b-n}.
\end{theorem}

The proof of the above theorem follows the ideas of Yao and Davis \cite{DY:1986} who showed that for i.i.d. standard normally distributed data the likelihood ratio converges in distribution to a Gumbel extreme value distribution. An important tool in the proof is the celebrated Darling-Erd\H{o}s theorem on the asymptotic distribution of $\max_{1\leq k\leq n} \frac{1}{\sqrt{k}} |\sum_{i=1}^k X_i|$. Theorem~4 below  establishes such a result for dependent data. Such theorems have been proved before, see e.g. Shorack \cite{S:1979}, but not under the conditions required in the present paper. 

\begin{theorem}
\label{thm:DE}
Let $(X_i)_{i\geq 1}$ be an $\alpha$-mixing strictly stationary process satisfying \eqref{eq:mixing-3} for some $2<p\leq 3$, and let $S_n:=\sum_{i=1}^n X_i$ denote the partial sum process.
 Then 
\begin{equation}
\frac{\sqrt{2\log\log n}}{\sigma}
\max_{1\leq k\leq n} \frac{|S_k|}{\sqrt{k}}-b_n \stackrel{\mathcal{D}}{\longrightarrow} G,
\label{eq:darling-erdos}
\end{equation}
where $G(x) =\exp(-\exp(-x))$ denotes the Gumbel extreme value distribution function, and where $\sigma$ and  $b_n$ is defined in \eqref{eq:lrv} and \eqref{eq:b-n}, respectively.
\end{theorem}
The proof of Theorem~4, presented in Section~5.4 below, follows the ideas of Shorack \cite{S:1979} who proved that an almost sure invariance principle together with a suitable maximal inequality implies the Darling-Erd\H{o}s theorem.

\section{Simulations}
\begin{figure}[h!]
\resizebox{\linewidth}{!}{\input{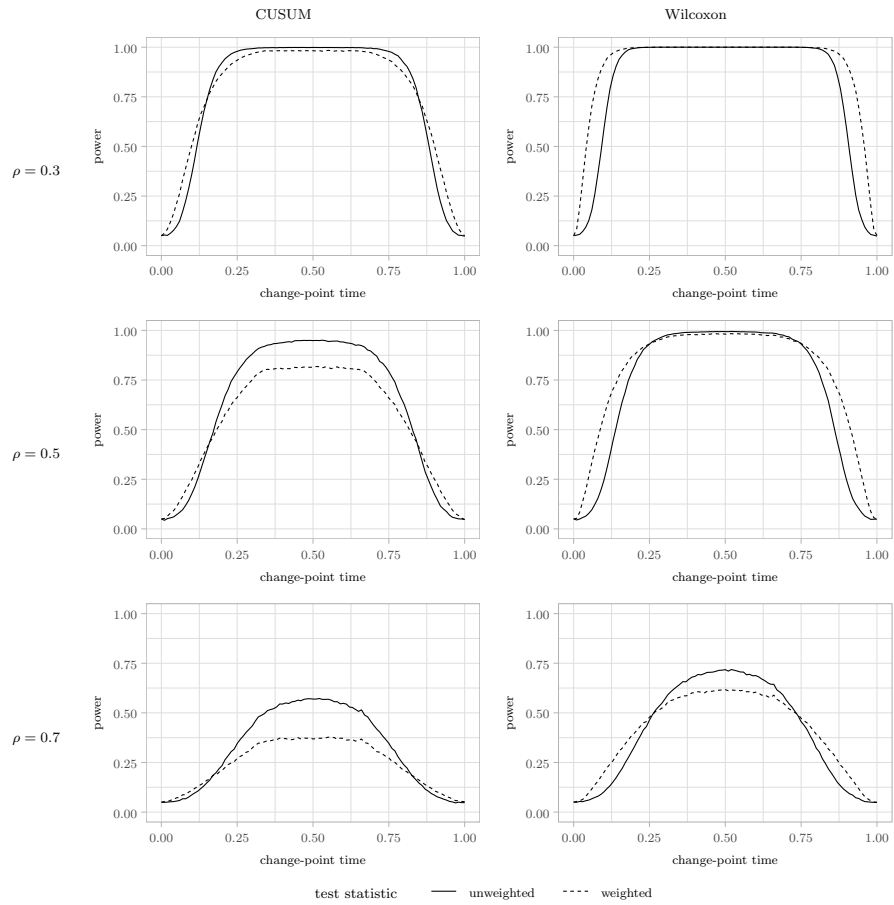}}
\captionof{figure}{Size-corrected power depending on the change-point time $\tau = [ \frac{k^*}{n} ].$ The data is generated from an AR(1) process with $t(3)$-distributed innovations and correlation coefficients $\rho=0.3,0.5,0.7$. The simulations are based on 20000 runs with sample size $n=800$ and shift height $\Delta=1.$ }
\label{fig.power.ar}
\end{figure}

In this section we present some simulation results for the weighted and the unweighted test statistic. We compare the power, the empirical size and the critical values, and consider the CUSUM and the Wilcoxon kernel, namely $h(x,y)=y-x$ and $h(x,y)=1_{\{ x<y\}}-\frac{1}{2}$.  
\begin{remark}
It is easy to see that the Wilcoxon kernel satisfies the variation condition. Let $h(x,y)=1_{\{x\leq y \}}-\frac{1}{2}$. Then, for $\varepsilon>0$, 
\begin{align*}
&\sup_{  \substack{\|(x,y)-(X,Y)\| \leq \varepsilon   }}    |h(x,y)-h(x',y')| \\
=&  \sup_{  \substack{\|(x,y)-(X,Y)\| \leq \varepsilon   }}    |1_{\{x\leq y\}}-1_{\{x'\leq y' \}}| \\
=& \begin{cases} 1 &~ \text{for}~X-Y \in (-\sqrt{2} \varepsilon, \sqrt{2} \varepsilon], \\ 0 &~\text{else}.  \end{cases}
\end{align*}
If the density of $X-Y$ is bounded, then 
\begin{align*}
& E \Big(  \sup_{  \substack{\|(x,y)-(X,Y)\| \leq \varepsilon   }}    |1_{\{x\leq y\}}-1_{\{x'\leq y' \}}| \Big)\\
\leq & P(X-Y \in(-\sqrt{2} \varepsilon, \sqrt{2} \varepsilon] ) \leq L \varepsilon.
\end{align*}
\end{remark}
Let $T_{n}^{\scriptscriptstyle{C}}$ denote the CUSUM and $T_n^{\scriptscriptstyle{WC}} $ the weighted CUSUM test statistic and let $T_{n}^{\scriptscriptstyle{W}}$ and $T_{n}^{\scriptscriptstyle{WW}}$ denote the Wilcoxon and the weighted Wilcoxon test statistic, all properly centered and normalized, i.e.\
\begin{align*}
T_{n}^{\scriptscriptstyle{C}}&=\frac{1}{ n^{3/2} \sigma_{\scriptscriptstyle{C}}}\max_{1\leq k \leq n-1}\left| \sum_{i=1}^k \sum_{j=k+1}^n (X_j-X_i) \right|, \\
{T}_{n}^{\scriptscriptstyle{WC}} &= \frac{\sqrt{2\log \log n}}{\sigma_{\scriptscriptstyle{C}}} \max_{1\leq k \leq n-1}  \frac{1}{\sqrt{k(n-k)n}} \left| \sum_{i=1}^k \sum_{j=k+1}^n  (X_j-X_i) \right|-b_n, \\
T_{n}^{\scriptscriptstyle{W}}&=\frac{1}{ n^{3/2} \sigma_{\scriptscriptstyle{W}}}\max_{1\leq k \leq n-1}\left| \sum_{i=1}^k \sum_{j=k+1}^n \left(1_{\{ X_i<X_j\}}-\frac{1}{2}\right) \right|, \\
{T}_{n}^{\scriptscriptstyle{WW}} &= \frac{\sqrt{2\log \log n}}{\sigma_{\scriptscriptstyle{W}}} \max_{1\leq k \leq n-1}  \frac{1}{\sqrt{k(n-k)n}} \left| \sum_{i=1}^k \sum_{j=k+1}^n  \left(1_{\{ X_i<X_j\}}-\frac{1}{2}\right) \right|-b_n.
\end{align*}
For $0< \alpha < 1$ we define $c_i(\alpha),~i\in \{1,2,3,4\},$ such that 
\begin{align*}
P\left( T_{n}^{\scriptscriptstyle{C}} > c_{1}( \alpha) \right)  =\alpha, \ \
 P\left({T}_{n}^{\scriptscriptstyle{WC}} > c_{2}(\alpha)\right) =\alpha, \\
P\left(T_{n}^{\scriptscriptstyle{W}} > c_{3}(\alpha)\right) =\alpha, \ \
P \left( {T}_{n}^{\scriptscriptstyle{WW}} > c_{4}(\alpha)\right) =\alpha.
\end{align*}
Most of the simulation study is based on i.i.d.\ standard normally distributed observations. In this case
\begin{equation*}
 \sigma_{\scriptscriptstyle{C}}^2=1 \ \text{and} \ \sigma_{\scriptscriptstyle{W}}^2=\frac{1}{12}.
\end{equation*}
In Figure~\ref{fig.power.ar} we we will consider dependent observations. In this case the long run variance has to be estimated. We use a subsampling estimator introduced by Dehling, Fried, Sharipov, Vogel and Wornowizki \cite{DFSVW:2013}. To achieve consistency under the alternative, we split the data into three disjoint sub-sequences of similar length and use the median of the resulting three separate estimations, see Dehling, Fried and Wendler \cite{DFW:2020}.

\begin{figure}[h!]
\resizebox{\linewidth}{!}{\input{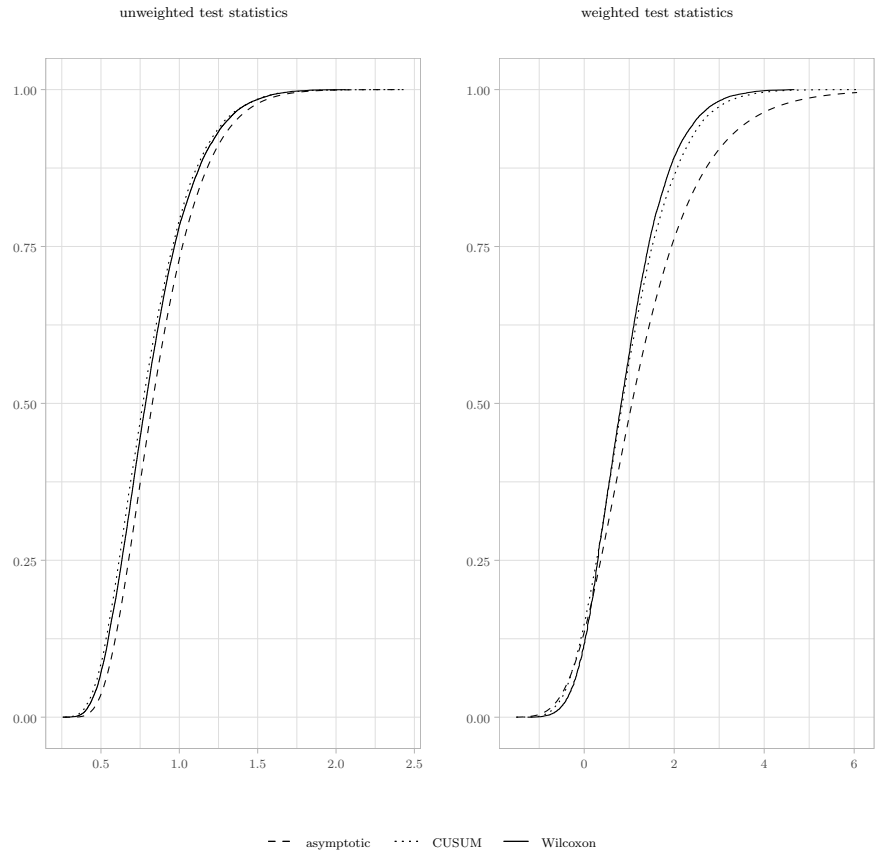}}
\captionof{figure}{Empirical distribution functions under the hypothesis for the unweighted and weighted CUSUM and Wilcoxon test statistics, compared to the distribution functions of the asymptotic distributions. The simulations are based on 100 i.i.d.\ standard normally distributed observations and 20000 runs.}
\label{fig.distr.fkt}
\end{figure}

First, we have simulated the critical values $c_i(\alpha)$ and compared them to the asymptotic ones. The results are summarized in Table~\ref{table.quantiles}. For the unweighted test statistics the simulated critical values are almost in agreement with the asymptotic ones, whereas for the weighted test statistics the difference is larger. An overview is also given in Figure~\ref{fig.distr.fkt}, which shows the different empirical distribution functions compared to the asymptotic ones. On the left hand side the empirical distribution function for the CUSUM and Wilcoxon test statistic is compared to the distribution function of the Kolmogorov-Smirnov distribution, and on the right hand side the distribution functions for the weighted CUSUM and weighted Wilcoxon test statistics are compared to the distribution function of the Gumbel distribution with location parameter $\log(2)$ and scale parameter 1. 

\begin{table}[h!]
\begin{tabular}{@{}rccccccccccc@{}}
\toprule
& & \multicolumn{5}{c}{$\alpha=0.05$}  & \multicolumn{5}{c}{$\alpha=0.10$} \\ \cmidrule(l){3-7} \cmidrule(l){8-12} 
    & n          & 100   & 200   & 400    & 800  & $\infty$ & 100   & 200    & 400   & 800 & $\infty$\\ \midrule 
 C & $c_1(\alpha)$ & 1.30 & 1.33 & 1.33  & 1.33  & 1.36  & 1.17 & 1.18  & 1.20 & 1.20 & 1.22 \\  
 WC & $c_2(\alpha)$ & 2.64 & 2.72 & 2.76  & 2.82  & 3.66  & 2.20 & 2.26  & 2.28 & 2.32 & 2.94 \\
 W & $c_3(\alpha)$ & 1.30 & 1.32  & 1.33 & 1.33  & 1.36  & 1.17 & 1.19   & 1.20 & 1.21 & 1.22 \\
 WW & $c_4(\alpha)$ & 2.47 & 2.56  & 2.60 & 2.65  & 3.66  & 2.05 & 2.13   & 2.18 & 2.18 & 2.94 \\
   
\bottomrule
\end{tabular}
\vspace{2mm}
\caption[Table]{Empirical and asymptotic critical values for the CUSUM \textup{(C)}, weighted CUSUM \textup{(WC)}, Wilcoxon \textup{(W)} and weighted Wilcoxon \textup{(WW)} test statistics. The empirical critical values are based on i.i.d.\ standard normally distributed observations and 20000 simulation runs.} 
\label{table.quantiles} 
\end{table}

   

Next, we evaluate the performance of the test statistics by computing the empirical sizes and the power. Table~\ref{emp.size.cusum.wilc} presents the empirical sizes for the unweighted and weighted CUSUM and Wilcoxon test statistics. The empirical sizes are lower than the nominal size in all cases considered. For the unweighted test statistics the size distortion shrinks to zero as the sample size increases, whereas for the weighted test statistics the difference between the empirical size and the nominal size is much larger for all considered sample sizes. 

\vspace{1em}
\begin{table}[h!]
\begin{tabular}{@{}cccccc@{}}
\toprule
 $n$ &    $\alpha$   &  $T_{n}^{\scriptscriptstyle{C}}$   & $T_{n}^{\scriptscriptstyle{WC}}$ & $T_{n}^{\scriptscriptstyle{W}}$ & ${T}_{n}^{\scriptscriptstyle{WW}}$  \\ \midrule 
 200      & 0.05     & 0.039   & 0.010         & 0.042   &  0.006  \\
 400      & 0.05     & 0.046   & 0.010         & 0.046   &  0.010\\
 800      & 0.05     & 0.049   & 0.014        & 0.047  &  0.012 \\
 1600     & 0.05     & 0.045   & 0.015         & 0.045   &  0.015 \\
 200	      & 0.10     & 0.085   & 0.035         & 0.088   &  0.023 \\
 400      & 0.10     & 0.085   & 0.039         & 0.085   &  0.039 \\
 800      & 0.10     & 0.095   & 0.043         & 0.091   &  0.034 \\
 1600     & 0.10     & 0.091   & 0.043         & 0.091   &  0.043 \\ 
\bottomrule
\end{tabular}
\vspace{2mm}
\caption[Table]{Empirical size based on $n$ i.i.d.\ standard normally distributed observations and 5000 simulation runs.} 
\label{emp.size.cusum.wilc} 
\end{table}

\begin{figure}[h!]
\resizebox{\linewidth}{!}{\input{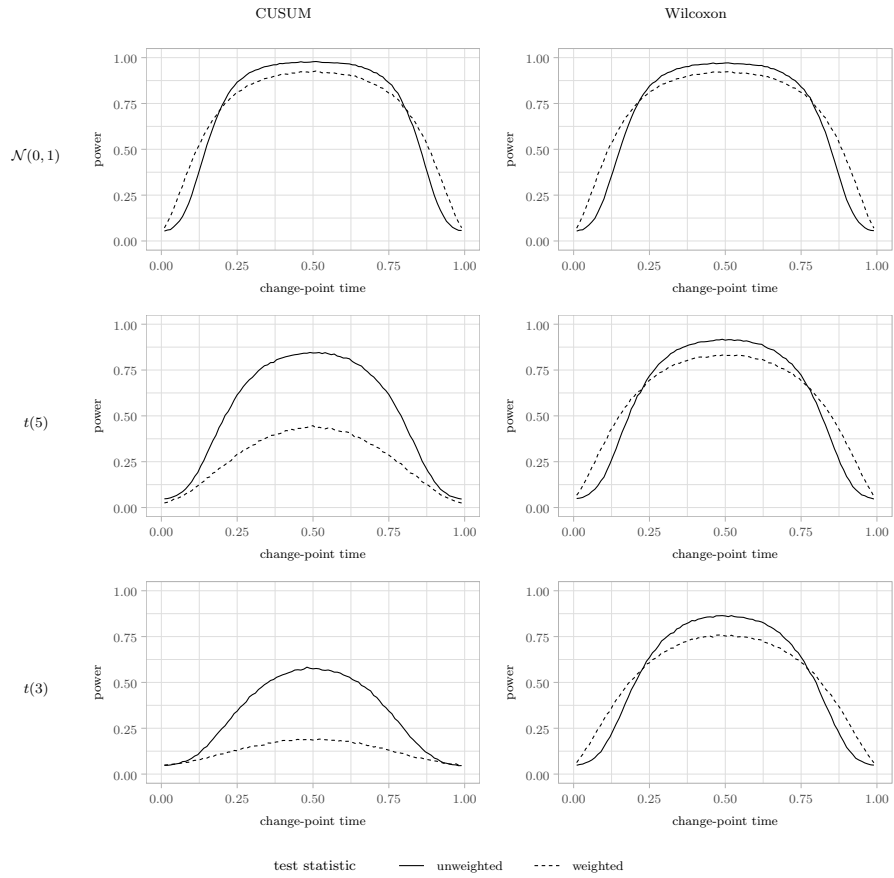}}
\captionof{figure}{Size-corrected power depending on the change-point time $\tau = [ \frac{k^*}{n} ].$ The simulations are based on 800 i.i.d.\ standard normally distributed observations with shift height $\Delta=0.3$ and 20000 simulation runs.} 
\label{fig.power}
\end{figure}

Figure~\ref{fig.power} shows the size-corrected power of the CUSUM, weighted CUSUM and the Wilcoxon and weighted Wilcoxon test statistic for all possible change-point times with shift height $\Delta=0.3.$ The simulations are based on i.i.d.\ standard normally distributed observations. It is clear that the power of all test statistics improves as the change-point is more central. In the case of early and late changes the weighted test statistics have better power than the unweighted ones, whereas they have less power for changes in the middle of the time period. For standard normally distributed data the power of the CUSUM and the Wilcoxon test statistics differs only marginally. Considering data with heavier tails, such as $t(5)$ or $t(3)$ distributed observations, one can see that the CUSUM test statistics lose more power, especially the weighted CUSUM. In both cases, the Wilcoxon test statistics have better power. 

Figure~\ref{fig.power.ar} shows analogous power plots, based on AR(1)-processes with $t(3)$-distributed innovations and correlation coefficients $\rho=0.3,0.5,0.7$. It is easy to see that a higher correlation results in a lower power. Again, in the case of early and late changes the weighted test statistics have better power than the unweighted ones. Compared to the CUSUM test statistics, the Wilcoxon test statistics always have greater power.  

\clearpage

\section{Data example}

\begin{figure}[h!]
\resizebox{\linewidth}{!}{\input{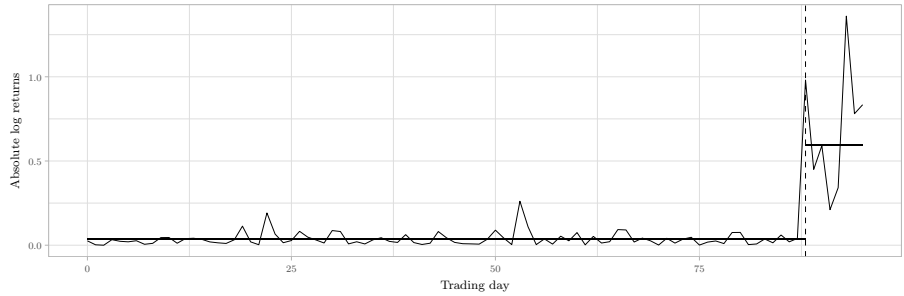}}
\captionof{figure}{Absolute log returns of the Wirecard stock price (Febraury 10, 2020 -- June 26, 2020). Change-point detected by the weighted Wilcoxon test (vertical dashed line), and average between the detected change-point (horizontal black lines).} 
\label{wirecard}
\end{figure}

As an application we analyze the daily absolute log returns of the closing Wirecard stock prices (currency in EUR, downloaded from \url{https://de.finance.yahoo.com/quote/WDI.DU/history?p=WDI.DU} on June 14, 2021). We consider the time period February 10, 2020, to June 26, 2020 , which is 19 weeks and 95 observations (trading time from monday to friday). The absolute log returns in this time period are visualized in Figure \ref{wirecard}. 
    
As the focus of this paper is on robust tests, we apply the weighted and non weighted Wilcoxon test to the data. The long run variance is estimated by the same subsampling estimator used in the simulation study. Considering a significance level of 5\%, the weighted Wilcoxon test statistic rejects the null hypothesis of a constant mean in the absolute log returns. In contrast, the unweighted Wilcoxon test doesn't detect any significant change (c.f. Table \ref{wirecardtable}). The weighted Wilcoxon test detects a change at observation 88, which correspond to June 18, 2020. On that day Wirecard reported that, following an audit by Ernst \& Young, about 1.9 billion euro was missing in certain trust accounts.    
\\

\begin{table}[h!]
\begin{tabular}{@{}cccccc@{}}
\toprule
 &   Wilcoxon &  weighted Wilcoxon    \\ \midrule 
 Test statistic      & 1.28     & 4.96   \\
 Critical value      & 1.36     & 3.66   \\
\bottomrule
\end{tabular}
\vspace{2mm}
\caption[Table]{Unweigted and weighted Wilcoxon test statistics together with the corresponding asymptotic critical values for the Wirecard stock price data, considering a significance level of 5\%.} 
\label{wirecardtable} 
\end{table}

\section{Proofs}

\subsection{Proof of Theorem~\ref{thm:mainresult}} \label{proof:mainresult}
First, we show in the following proposition that the degenerate part is asymptotically negligible.
\begin{proposition}
\label{prop.deg.part}
Let $\Psi$ be the kernel given by Hoeffding's decomposition of $h$ in \eqref{hoeffding}, satisfying the variation condition. Moreover, we assume that \eqref{eq:mix-2} holds. Then, under the null hypothesis $H_0$, as $n\rightarrow \infty$
\begin{align*}
\max_{1\leq k \leq n-1} \sqrt{\frac{\log\log n}{k(n-k)n}} \left|  \sum_{i=1}^k  \sum_{j=k+1}^n \Psi (X_i,X_j) \right| \cpr 0. 
\end{align*} 
\end{proposition}
\begin{proof}
We can split the maximum into the stretch up to $\sqrt{n}$, the stretch between $\sqrt{n}$ and $n-\sqrt{n}$ and the stretch after $n-\sqrt{n}$, such that
\begin{align} 
& \max_{1\leq k \leq n-1} \sqrt{\frac{\log\log n}{k(n-k)n}} \left| \sum_{i=1}^k \sum_{j=k+1}^n \Psi (X_i,X_j) \right| \notag \\
\leq 2 & \max_{\sqrt{n}\leq k<n-\sqrt{n}} \sqrt{\frac{\log\log n}{n^{5/2}}} \left| \sum_{i=1}^k \sum_{j=k+1}^n \Psi (X_i,X_j)\right| \label{max1}\\
+ & \max_{k<\sqrt{n}}  \sqrt{\frac{\log\log n}{n^2}} \left| \sum_{i=1}^k \sum_{j=k+1}^n \Psi (X_i,X_j) \right| \label{max2} \\
+ & \max_{k\geq n-\sqrt{n}} \sqrt{\frac{\log\log n}{n^2}} \left| \sum_{i=1}^k \sum_{j=k+1}^n \Psi (X_i,X_j) \right|. \label{max3}
\end{align}
Now we can deal with every single maximum. To show that these maxima converge in probability to $0$, we use Theorem A of Serfling \cite{S:1970}. To apply that theorem, we need a functional $g(F_{a,n})$ depending on the joint distribution of a vector $(Y_{a+1},\dots,Y_{a+n})$ of $n$ random variables, and satisfying 
\[g(F_{a,k})+g(F_{a+k,l}) \leq g(F_{a,k+l}),~1\leq k<k+l, \]
such that 
\[ E|S_{a,l}|^2 \leq g(F_{a,l}),~\text{for all} ~l\geq 1,\]
where 
\[S_{a,l}=\sum_{i=a+1}^{a+l} Y_i.\] 
Then it follows
\[ E\left( \max_{1 \leq i \leq l } |S_{a,i}| \right)^2 \leq\left( \frac{\log(2l)}{\log(2)} \right)^2 \left( g(F_{a,l}) \right)^2. \] 
We set $V_k=\sum_{i=1}^{k} \sum_{j=k+1}^n \Psi (X_i,X_j).$ By Lemma~\ref{Lemma_DFGW}
\begin{equation*}
E(V_m-V_k)^2\leq C_1 |m-k|n,
\end{equation*}
where $C_1$ is a constant. Define $Y_i:=V_i-V_{i-1},~i=1,\dots,n$ and set $V_0=0$. Then one has $V_m=\sum_{i=1}^{m} Y_i$ and $V_m-V_k=\sum_{i=k+1}^{m} Y_i$. As $k\leq m$ we can set $m:=k+l,~l\in \mathbb{N}.$ Thus 
\[E\left(\sum_{i=k+1}^{k+l} Y_i\right)^2\leq C_1ln. \] 
With
\[g(F_{a,l}) =C_1ln,\]
where $C_1$ is a constant, the required conditions are satisfied. We have
\[g(F_{a,k})+g(F_{a+k,l})=C_1kn + C_1ln= C_1n(k+l)=g(F_{a,k+l})\]
and
\[E|S_{a,l}|^2=E \left( \sum_{i=a+1}^{a+l} Y_i \right) \leq C_1ln=g(F_{a,l}).\]
Thus
\[E \left( \max_{1\leq j \leq l} \left| \sum_{i=k+1}^{k+j} Y_i \right|\right) ^2 \leq \left( \frac{\log (2l)}{\log(2)}\right)^2C_1ln. \] 
For $k=0$ and $l=\sqrt{n}$ we get 
\begin{equation*}
E\left(\max_{1\leq j\leq \sqrt{n}} |V_j|\right)^2= E \left( \max_{1\leq j \leq \sqrt{n}} \left| \sum_{i=1}^{j} Y_i \right|\right)^2 \leq \left( \frac{\log \left(2\sqrt{n}\right)}{\log(2)}\right)^2C_1 n^{3/2}.
\end{equation*}
Thus we have 
\begin{equation*}
E\left(\max_{1\leq k\leq \sqrt{n}} \sqrt{\frac{\log \log n}{n^2}}|V_k|\right)^2\leq \frac{\log \log n}{\sqrt{n}} \left( \frac{\log \left(2\sqrt{n}\right)}{\log(2)}\right)^2C_1 .
\end{equation*}
This goes to 0 for $n\rightarrow \infty.$ Applying Chebyshev's inequality, we obtain that \eqref{max2} converges to $0$ in probability as $n\rightarrow \infty$. By stationarity, this also holds for \eqref{max3}. An analogous procedure leads to
\begin{equation*}
E\left(\max_{1\leq k\leq n-\sqrt{n}} \sqrt{\frac{\log \log n}{n^{5/2}}}|V_k|\right)^2\leq \frac{(n-\sqrt{n})\log \log n}{n^{3/2}} \left( \frac{\log \left(2\sqrt{n}\right)}{\log(2)}\right)^2C_1 .
\end{equation*}
As this goes to zero for $n\rightarrow \infty$ and as
\begin{equation*}
\max_{\sqrt{n} \leq k\leq n-\sqrt{n}} \sqrt{\frac{\log \log n}{n^{5/2}}}|V_k| \leq \max_{1\leq k\leq n-\sqrt{n}} \sqrt{\frac{\log \log n}{n^{5/2}}}|V_k|,
\end{equation*}
we obtain that \eqref{max1} converges in probability to 0.  

\end{proof}

\begin{proof}[Proof of Theorem~\ref{thm:mainresult}]

In the same way as in \eqref{decomp}, we apply Hoeffding's decomposition to the kernel $h(x,y)$ and get
\begin{align*}
 &\frac{1}{\sqrt{k(n-k)n}}  \sum_{i=1}^k \sum_{j=k+1}^n h(X_i,X_j) \\ 
=& \sqrt{\frac{n}{k(n-k)}}\left( \sum_{i=1}^k h_1(X_i)-\frac{k}{n} \sum_{i=1}^n h_1(X_i) \right)+ \frac{1}{\sqrt{k(n-k)n}}\sum_{i=1}^k \sum_{j=k+1}^n \Psi(X_i,X_j). 
\end{align*}

As the variation condition holds for the kernel $h$, it also holds for $\Psi.$ In order to be able to apply Proposition \ref{prop.deg.part}, we need to verify assumption \eqref{eq:mix-2}. From \eqref{ass.epsilon}, we obtain
\begin{align*}
\int_0^1( \alpha^{-1}(u))^{4+\varepsilon} \, Q_{|X|}(u) du<\infty,
\end{align*}
which is equivalent to 
\begin{align*}
\sum_{k=1}^{\infty} k^{3+ \varepsilon} \int_{0}^{\alpha(k)} Q_{|X|}(u)du < \infty;
\end{align*}
see Remark 2.2 in Merlev\`ede and Rio \cite{MR:2012} or Annex C in Rio \cite{Rio:2017}. Now, we define
\begin{align*}
q_k:=\int_{0}^{\alpha(k)} Q_{|X|}(u)du.
\end{align*}
As $q_k$ is nonincreasing and by assumption, we obtain 
\begin{align*}
n^{4+\varepsilon} q_n  \leq c_1 \sum_{k=1}^n k^{3+\varepsilon} q_n  
  \leq c_1 \sum_{k=1}^n k^{3+\varepsilon} q_k 
  \leq c_1 \sum_{k=1}^{\infty} k^{3+\varepsilon} q_k < \infty,
\end{align*}
where $c_1$ is a constant. Hence $n^{\frac{(4+\varepsilon)\delta}{3+2\delta}} q_n^{\frac{\delta}{3+2\delta}} $ is bounded by some constant $c_2$, and thus
\begin{align*}
n q_n^{\frac{\delta}{3+2\delta}} \leq c_2 \frac{1}{n^{1+\frac{\varepsilon \delta -6}{3+2 \delta}}},
\end{align*}
which implies
\begin{align*}
\sum_{k=1}^{\infty} k q_k^{\frac{\delta}{3+2\delta}} < \infty.
\end{align*}
Thus, all conditions of Proposition~\ref{prop.deg.part} are satisfied and we can conclude that
\begin{align*}
\max_{1\leq k \leq n-1} \sqrt{\frac{\log \log n}{k(n-k)n}} \left| \sum_{i=1}^k \sum_{j=k+1}^n \Psi(X_i,X_j) \right| \cpr 0.
\end{align*}
Hence, by Slutsky's theorem, it remains to show that
\begin{equation*}
\max_{1\leq k \leq n-1} \sqrt{\frac{n \log \log n}{k(n-k)}} \left|  \sum_{i=1}^k h_1(X_i)-\frac{k}{n} \sum_{i=1}^n h_1(X_i) \right|- b_n 
\end{equation*}
converges in distribution to the desired extreme value distribution. This follows from Theorem~\ref{thm:tied-down} with $S_k=\sum_{i=1}^k h_1(X_i)$.
\end{proof}

\subsection{Proof of Theorem~\ref{thm:consistency}} \label{proof:consistency}

Without loss of generality, we assume that $\Delta>0$. Since
\begin{align*}
\max_{1\leq k \leq n-1} \! \frac{1}{\sqrt{k(n-k)n}} \! \Big| \! \sum_{i=1}^k \! \sum_{j=k+1}^n \! h(X_i,X_j) \! \Big| \! \geq \! \frac{1}{\sqrt{k^*(n-k^*)n}} \! \Big| \! \sum_{i=1}^{k^*} \! \sum_{j=k^*+1}^n \! h(X_i,X_j) \! \Big|,
\end{align*}
it suffices to show that 
\begin{align*}
\frac{1}{\sqrt{(\log \log n) k^*(n-k^*)n}}  \left| \sum_{i=1}^{k^*} \sum_{j=k^*+1}^n h(X_i,X_j) \right| \cpr \infty.
\end{align*}
We obtain the following Hoeffding decomposition
\[
 h(X_i,X_j)=\Delta +h_1(X_i) +h_2(X_j) +\Psi(X_i,X_j), \quad 1\leq i \leq k^\ast<j\leq n,
\]
where $\Delta$ is given as in Theorem~\ref{thm:consistency}, and where
\begin{align*}
 h_1(x)&:= Eh(x,X_n)-\Delta, \\
 h_2(y)&:= Eh(X_1,y) -\Delta.
\end{align*}
It holds 
\begin{align*}
\frac{T_n}{\sqrt{\log \log n}} \! \geq \! &  \frac{1}{\sqrt{(\log \log n)k^*(n-k^*)n}}  \!  \sum_{i=1}^{k^*} \! \sum_{j=k^*+1}^{n} \! (\Delta \!+\! h_{1}(X_i) \! + \! h_{2}(X_j) \! + \! \Psi(X_i,X_j) )  \\
   = & \sqrt{\frac{k^*(n-k^*)}{(\log \log n) n}} \Delta+ \sqrt{\frac{n-k^*}{(\log \log n) k^*n}} \sum_{i=1}^{k^*} h_{1}(X_i) \\
   &+ \sqrt{\frac{k^*}{(\log \log n)(n-k^*)n}} \sum_{j=k^*+1}^{n} h_{2}(X_j) \\
   &+ \frac{1}{\sqrt{(\log \log n) k^*(n-k^*)n}} \sum_{i=1}^{k^*} \sum_{j=k^*+1}^{n} \Psi(X_i,X_j).
\end{align*}
By the law of iterated logarithm  
\begin{align*}
\sqrt{\frac{n-k^*}{ k^*n \log \log n}} \Big| \sum_{i=1}^{k^*} h_{1}(X_i) \Big| \leq & \frac{1}{\sqrt{k^* \log \log n}} \Big| \sum_{i=1}^{k^*} h_{1}(X_i) \Big| =O(1)~a.s., \\
\sqrt{\frac{k^*}{(n-k^*)n \log \log n}} \Big| \sum_{j=k^*+1}^{n} h_{2}(X_j)\Big| \leq  &\frac{1}{\sqrt{(n-k^* \log \log n)}} \Big| \sum_{j=k^*+1}^{n} h_{2}(X_j)\Big| \\\overset{\mathcal{D}}{=}&\frac{1}{\sqrt{(n-k^*)\log \log n}}\Big| \sum_{j=1}^{n-k^*} h_{2}(X_{j+k^*}) \Big| \\ =&  O_p(1).
\end{align*}
It remains to show that
\begin{align*}
\frac{1}{\sqrt{(\log \log n) k^*(n-k^*)n}} \sum_{i=1}^{k^*} \sum_{j=k^*+1}^{n} \Psi(X_i,X_j) \cpr 0,
\end{align*}
as $n\rightarrow \infty.$ By Lemma~\ref{Lemma_DFGW} we have 
\begin{align*}
E\left( \sum_{i=1}^{k^*} \sum_{j=k^*+1}^{n} \Psi(X_i,X_j) \right)^2 \leq C(n-k^*)k^*,
\end{align*}
for some constant $C$. We get 
\begin{align*}
E\left(\frac{1}{\sqrt{(\log \log n) k^*(n-k^*)n}} \sum_{i=1}^{k^*} \sum_{j=k^*+1}^{n} \Psi(X_i,X_j) \right)^2 \leq \frac{C}{(\log \log n)n} \rightarrow 0, 
\end{align*}
as $n\rightarrow \infty,$ which completes the proof.

\subsection{Proof of Theorem~\ref{thm:tied-down}}
\label{proof:tied-down}

For the proof of Theorem~\ref{thm:tied-down}, we need the following lemmas, which all hold under the assumptions of Theorem~\ref{thm:tied-down}.

\begin{lemma}
\label{lemma2.1}
\begin{align*}
\max_{1\leq k \leq  \left[\frac{n}{\log n}\right]} \sqrt{\frac{n}{k(n-k)}} \left| S_k-\frac{k}{n}S_n \right| -\max_{1\leq k \leq  \left[\frac{n}{\log n}\right]} \left| \frac{S_k}{\sqrt{k}} \right| =o_p\left( \frac{1}{\sqrt{2\log \log n}} \right)
\end{align*}
\end{lemma}
\begin{proof}
For $n$ large, we have 
\begin{align}
& \left| \max_{1\leq k \leq  \left[\frac{n}{\log n}\right]} \sqrt{\frac{n}{k(n-k)}} \left| S_k-\frac{k}{n}S_n \right| -\max_{1\leq k \leq  \left[\frac{n}{\log n}\right]} \left| \frac{S_k}{\sqrt{k}} \right| \right| \notag \\
\leq & \max_{1\leq k \leq  \left[\frac{n}{\log n}\right]} \left|  \sqrt{\frac{n}{k(n-k)}} \left| S_k-\frac{k}{n}S_n \right| - \left| \frac{S_k}{\sqrt{k}} \right| \right| \notag \\
\leq & \max_{1\leq k \leq  \left[\frac{n}{\log n}\right]} \left|  \sqrt{\frac{n}{k(n-k)}}  \left( S_k-\frac{k}{n}S_n \right)  -  \frac{S_k}{\sqrt{k}}  \right| \notag \\
 = & \max_{1\leq k \leq  \left[\frac{n}{\log n}\right]} \left| \left( \sqrt{\frac{n}{(n-k)}}-1 \right) \frac{S_k}{\sqrt{k}}  - \sqrt{\frac{k}{n-k}} \frac{S_n}{\sqrt{n}}  \right| \notag \\
 \leq & \max_{1\leq k \leq  \left[\frac{n}{\log n}\right]} \left(\left| \left( \sqrt{\frac{n}{(n-k)}}-1 \right) \frac{S_k}{\sqrt{k}} \right| + \left|\sqrt{\frac{k}{n-k}} \frac{S_n}{\sqrt{n}} \right| \right) \label{ineq}
\end{align}
For large $n$ and $1\leq k \leq \left[\frac{n}{\log n}\right]$ the inequalities 
\begin{align*}
\sqrt{\frac{n}{n-k}}-1 \leq \frac{k}{n}~\text{and}~\sqrt{\frac{k}{n-k}}\leq 2\sqrt{\frac{k}{n}}
\end{align*}
are satisfied. So we have that \eqref{ineq} is less or equal to  
\begin{align*}
  & \max_{1\leq k \leq  \left[\frac{n}{\log n}\right]} \left(\left| \frac{k}{n}    \frac{S_k}{\sqrt{k}} \right| + \left|2 \sqrt{\frac{k}{n}} \frac{S_n}{\sqrt{n}} \right| \right) \\
= & \max_{1\leq k \leq  \left[\frac{n}{\log n}\right]} \sqrt{\frac{k}{n}} \left(\left|   \frac{S_k}{\sqrt{n}} \right| + 2\left| \frac{S_n}{\sqrt{n}} \right| \right) \\
\leq & \frac{1}{\sqrt {\log n}}\max_{1\leq k \leq  \left[\frac{n}{\log n}\right]}  \left(\left|   \frac{S_k}{\sqrt{n}} \right| + 2\left| \frac{S_n}{\sqrt{n}} \right| \right) \\
= & O_p\left( \frac{1}{\sqrt {\log n}}\right),
\end{align*} 
where the last equality holds, as 
\begin{align*}
\max_{1\leq k \leq  \left[\frac{n}{\log n}\right]}  \left(\left|   \frac{S_k}{\sqrt{n}} \right| + 2\left| \frac{S_n}{\sqrt{n}} \right| \right)= O_p(1),
\end{align*}
 since $\max_{1\leq k \leq  \left[\frac{n}{\log n}\right]}  \left|   \frac{S_k}{\sqrt{n}} \right|$ and $ \left|  \frac{S_n}{\sqrt{n}} \right|$ converge in distribution. Now the claim follows, as $ ( \frac{1}{\sqrt{\log n}} ) /  ( \frac{1}{\sqrt{2\log \log n}} ) \rightarrow 0,$ for $n\rightarrow \infty.$ 
\end{proof}

\begin{lemma}
\label{lemma2.2}
\begin{align*}
\max_{\left[\frac{n}{\log n}\right] \leq k \leq  \left[\frac{n}{2}\right]} \sqrt{\frac{n}{k(n-k)}} \left| S_k-\frac{k}{n}S_n \right|= O_p\left(\sqrt{\log \log \log n } \right)
\end{align*}
\end{lemma}
\begin{proof} Since
\begin{align*}
& \max_{\left[\frac{n}{\log n}\right] \leq k \leq  \left[\frac{n}{2}\right]} \sqrt{\frac{n}{k(n-k)}} \left| S_k-\frac{k}{n}S_n \right| \\
\leq &  \max_{\left[\frac{n}{\log n}\right] \leq k \leq  \left[\frac{n}{2}\right]}2 \left|\frac{S_k}{\sqrt{k}} \right|+ \max_{\left[\frac{n}{\log n}\right] \leq k \leq  \left[\frac{n}{2}\right]} \left| \frac{S_n}{\sqrt{n}} \right|
\end{align*}
and the second summand is $O_p(1),$ it suffices to show that 
\begin{align*}
\max_{\left[\frac{n}{\log n}\right] \leq k \leq  \left[\frac{n}{2}\right]} \left|\frac{S_k}{\sqrt{k}} \right| = O_p\left( \sqrt{\log \log \log n } \right) .
\end{align*}
 By the almost sure invariance principle of Merlev\`ede and Rio (2012), one can find a Brownian motion $(W_t)_{t\geq 0}$ such that 
\begin{align*}
\frac{|S_t-W_t|}{\sqrt{t}}= O\left(t^{-\lambda}\right)~a.s.,
\end{align*}
for some $\lambda>0$. Hence
\begin{align*}
\max_{\left[\frac{n}{\log n}\right] \leq k \leq  \left[\frac{n}{2}\right]} \left|\frac{S_k-W_k}{\sqrt{k}} \right| =O\left(\left(\frac{\log n}{n} \right)^{\lambda} \right) ~ a.s.,
\end{align*}
and thus
\begin{align}
\label{221}
\max_{\left[\frac{n}{\log n}\right] \leq k \leq  \left[\frac{n}{2}\right]} \left|\frac{S_k}{\sqrt{k}} \right| \leq  \max_{\left[\frac{n}{\log n}\right] \leq k \leq  \left[\frac{n}{2}\right]} \left|\frac{W_k}{\sqrt{k}} \right| + \max_{\left[\frac{n}{\log n}\right] \leq k \leq  \left[\frac{n}{2}\right]} \left|\frac{S_k-W_k}{\sqrt{k}} \right|. 
\end{align}
Set $k=nr,~ r\in (0,1)$ and consider the first summand of the right-hand side of \eqref{221}. We have 
\begin{align*}
 \max_{\left[\frac{n}{\log n}\right] \leq nr \leq  \left[\frac{n}{2}\right]} \left|\frac{W_{nr}}{\sqrt{nr}} \right| 
\stackrel{\mathcal{D}}{=} \max_{\left[\frac{n}{\log n}\right] \leq nr \leq  \left[\frac{n}{2}\right]} \left|\frac{W_{r}}{\sqrt{r}} \right| \leq \max_{  \frac{1}{\log n} \leq r \leq  \frac{1}{2}} \left|\frac{W_{r}}{\sqrt{r}} \right| .
\end{align*}
By the law of the iterated logarithm, it holds for $r\rightarrow 0$
\begin{align*}
|W_r|=O_p\left(\sqrt{2r\log \log r^{-1}}\right).
\end{align*}
It follows that 
\begin{align*}
\max_{  \frac{1}{\log n} \leq r \leq  \frac{1}{2}} \left|\frac{W_{r}}{\sqrt{r}} \right| =O_p\left( \sqrt{\log \log \log n } \right),
\end{align*}
which completes the proof.
\end{proof}
For abbreviation, we define $a_n\!:=\! \sqrt{2 \log \log n}$. Recall the definition of $b_n$ \eqref{eq:b-n}.

\begin{lemma}
\label{lemma2.3}
\begin{align*}
\lim_{n\rightarrow \infty } P \left( \frac{a_n}{\sigma}  \max_{1 \leq k \leq  \left[\frac{n}{\log n}\right]} \left| \frac{S_k}{\sqrt{k}}\right|  -b_n \leq x \right) = \exp (-\exp (-x))
\end{align*}
\end{lemma}
\begin{proof}
By Theorem~\ref{thm:DE}, one has
\begin{align*}
\frac{a_{\left[\frac{n}{\log n}\right]}}{\sigma} \max_{1\leq k \leq \left[\frac{n}{\log n}\right] } \frac{|S_k|}{\sqrt{k}} - b_{\left[\frac{n}{\log n}\right]}  \stackrel{\mathcal{D}}{\longrightarrow} G.
\end{align*}
Since 
\begin{align*}
\frac{a_n}{a_{\left[\frac{n}{\log n}\right]}} \rightarrow 1, \ \ b_{\left[\frac{n}{\log n}\right]}\frac{ a_n}{a_{\left[\frac{n}{\log n}\right]}}- b_n \rightarrow 0,
\end{align*}
as $n\rightarrow \infty,$ we obtain with Slutsky's Theorem 
\begin{align*}
& \frac{a_n}{\sigma} \max_{1\leq k \leq \left[\frac{n}{\log n}\right] } \frac{|S_k|}{\sqrt{k}} -b_n  \\
= \  &  \frac{a_n}{a_{\left[\frac{n}{\log n}\right]}} \left( \frac{a_{\left[\frac{n}{\log n}\right]}}{\sigma} \max_{1\leq k \leq \left[\frac{n}{\log n}\right] } \frac{|S_k|}{\sqrt{k}} - b_{\left[\frac{n}{\log n}\right]} \right) + b_{\left[\frac{n}{\log n}\right]}\frac{ a_n}{a_{\left[\frac{n}{\log n}\right]}}- b_n  \claw G.
\end{align*}
\end{proof}

\begin{lemma}
\label{lemma2.4}
\begin{align*}
\max_{1 \leq k \leq  \left[\frac{n}{2}\right]} \sqrt{\frac{n}{k(n-k)}} \left| S_k-\frac{k}{n}S_n \right| - \max_{1\leq k \leq  \left[\frac{n}{\log n}\right]} \left| \frac{S_k}{\sqrt{k}} \right| =o_p\left( \frac{1}{\sqrt{2\log \log n}} \right)
\end{align*}
\end{lemma}
\begin{proof}
Since $\frac{1}{a_n^2}\rightarrow 0, $ we can conclude from Lemma~\ref{lemma2.3} and Slutsky's Theorem
\begin{align*}
\frac{1}{a_n^2}\left( a_n \max_{1\leq k \leq \left[\frac{n}{\log n}\right] } \frac{|S_k|}{\sqrt{k}} -b_n \right) \cpr 0.
\end{align*}
As $\frac{b_n}{a_n^2}\rightarrow 1,$ we get 
\begin{align*}
\frac{1}{a_n} \max_{1\leq k \leq \left[\frac{n}{\log n}\right] } \frac{|S_k|}{\sqrt{k}} \cpr 1.
\end{align*}
By Lemma~\ref{lemma2.1}, this implies 
\begin{align*}
\frac{1}{a_n} \max_{1\leq k \leq \left[\frac{n}{\log n}\right] } \sqrt{\frac{n}{k(n-k)}} \left| S_k-\frac{k}{n}S_n \right|  \cpr 1.
\end{align*}
By Lemma~\ref{lemma2.2}, we have
\begin{align*}
\frac{1}{a_n} \max_{ \left[\frac{n}{\log n}\right] \leq k \leq \left[\frac{n}{2}\right] } \sqrt{\frac{n}{k(n-k)}} \left| S_k-\frac{k}{n}S_n \right|  \cpr 0.
\end{align*}
Thus
\begin{align*}
&P\left( \max_{\left[\frac{n}{\log n}\right] \leq k \leq  \left[\frac{n}{2}\right]} \sqrt{\frac{n}{k(n-k)}} \left| S_k-\frac{k}{n}S_n \right| \geq \max_{1 \leq k \leq \left[\frac{n}{\log n}\right] } \sqrt{\frac{n}{k(n-k)}} \left| S_k-\frac{k}{n}S_n \right|  \right)\\ 
 &\longrightarrow 0,
\end{align*}
as $n\rightarrow \infty,$ and hence the lemma follows from Lemma~\ref{lemma2.1}.


\end{proof}

\begin{lemma}
\label{lemma2.5}
\begin{align*}
\max_{1 \leq n-k \leq  \left[\frac{n}{2}\right]} \! \sqrt{\frac{n}{k(n-k)}} \! \left|  \!S_k \! -\! \frac{k}{n}S_n \! \right|\! -\! \max_{1\leq n-k \leq \left[\frac{n}{\log n}\right]} \! \left| \! \frac{S_n \! - \!S_k}{\sqrt{n-k}} \! \right| =o_p\! \left(\! \frac{1}{\sqrt{2\log \log n}}\! \right)
\end{align*}
\end{lemma}
\begin{proof}
Follows from Lemma~\ref{lemma2.4} and the stationarity under the hypothesis. 
\end{proof}

Applying the above lemmas, we can now prove Theorem~\ref{thm:tied-down}.
\begin{proof}[Proof of Theorem~\ref{thm:tied-down}]
We have 
\begin{align*}
& P\left( \frac{a_n}{\sigma}  \max_{1\leq k \leq n-1} \sqrt{\frac{n}{k(n-k)}} \left| S_k-\frac{k}{n}S_n \right| - b_n \leq x \right) \\
=\ & P\left( \max_{1\leq k \leq \left[\frac{n}{2}\right]} \sqrt{\frac{n}{k(n-k)}} \left| S_k-\frac{k}{n}S_n \right|   \leq \frac{(x+b_n) \sigma}{a_n}, \right. \\ 
& ~~~~~~ \left. \max_{1\leq n-k \leq \left[\frac{n}{2}\right]} \sqrt{\frac{n}{k(n-k)}} \left| S_k-\frac{k}{n}S_n \right|  \leq \frac{(x+b_n)\sigma}{a_n}  \right).
\end{align*} 
From Lemma~\ref{lemma2.4} and Lemma~\ref{lemma2.5} we get 
\begin{align*}
\frac{a_n}{\sigma}  \max_{1\leq k \leq \left[\frac{n}{2}\right]} \sqrt{\frac{n}{k(n-k)}} \left| S_k-\frac{k}{n}S_n \right|- \frac{a_n}{\sigma} \max_{1\leq k \leq \left[\frac{n}{\log n}\right]}  \left| \frac{S_k}{\sqrt{k}} \right| \cpr 0
\end{align*}
and 
\begin{align*}
\frac{a_n}{\sigma}  \max_{1\leq n-k \leq \left[\frac{n}{2}\right]} \sqrt{\frac{n}{k(n-k)}} \left| S_k-\frac{k}{n}S_n \right|- \frac{a_n}{\sigma} \max_{1\leq n-k \leq \left[\frac{n}{\log n}\right]} \left| \frac{S_n-S_k}{\sqrt{n-k}} \right| \cpr 0.
\end{align*} Applying Lemma~\ref{lemma2.3} we get
\begin{align*}
\lim_{n \rightarrow \infty} P\left( \max_{1\leq k < \left[\frac{n}{\log n}\right]}  \left| \frac{S_k}{\sqrt{k}} \right|  \leq \frac{(x+b_n)\sigma}{a_n} \right) = \exp (-\exp(-x)).
\end{align*}
And again applying Lemma~\ref{lemma2.3} to the reversed time series, we get
\begin{align*}
\lim_{n \rightarrow \infty} P\left( \max_{1\leq n-k \leq \left[\frac{n}{\log n}\right]} \left| \frac{S_n-S_k}{\sqrt{n-k}} \right|  \leq \frac{(x+b_n)\sigma}{a_n}  \right)= \exp (-\exp(-x)).
\end{align*}
Due to the underlying $\alpha-$mixing process it holds
\begin{align*}
&\left| P\left( \max_{1\leq k \leq \left[\frac{n}{\log n}\right]}  \left| \frac{S_k}{\sqrt{k}} \right|  \leq \frac{(x+b_n)\sigma}{a_n},  \max_{1\leq n-k \leq \left[\frac{n}{\log n}\right]} \left| \frac{S_n-S_k}{\sqrt{n-k}} \right|  \leq \frac{(x+b_n)\sigma}{a_n}  \right) \right. \\ 
& \left. - P \! \left(\! \max_{ \! 1\leq k \leq \left[\frac{n}{\log n}\right]} \! \left|\! \frac{S_k}{\sqrt{k}}\! \right| \! \leq\! \frac{(\!x\!+\!b_n\!)\sigma}{a_n} \! \right)\! P\! \left(\! \max_{1\leq n-k \! \leq \left[\frac{n}{\log n}\right]}\! \left|\! \frac{S_n-S_k}{\sqrt{n-k}}\! \right| \! \leq\! \frac{(\!x\!+\!b_n\!)\sigma}{a_n} \! \right)\! \right|  \longrightarrow 0.
\end{align*}
All in all we get the desired result
\begin{align*}
\lim_{n \rightarrow \infty} P\left( \frac{a_n}{\sigma}  \max_{1\leq k \leq n-1} \sqrt{\frac{n}{k(n-k)}} \left| S_k-\frac{k}{n}S_n \right| - b_n \leq x \right)= \exp (-2\exp(-x)).
\end{align*}

\end{proof}

\subsection{Proof of Theorem~\ref{thm:DE}} \label{proof:DE}
First, we state two lemmas, both valid under the assumptions of Theorem~\ref{thm:DE}. We define the random variables $Y_n$ and $Z_n$ by
\begin{eqnarray*}
Y_n&=& \max_{1\leq k \leq (\log\log n)^\gamma} \frac{|S_k|}{\sqrt{k}} ,     \\
Z_n&=& \max_{(\log\log n)^\gamma\leq k \leq n} \frac{|S_k|}{\sqrt{k}},  
\end{eqnarray*}
for some $\gamma>\frac{1}{2}$. 
\begin{lemma}
\label{lem:y-n}
\begin{align*}
\sqrt{2\log\log n}\, Y_n -2\log\log n \cpr -\infty
\end{align*}
\end{lemma}
\begin{proof}
We apply the following maximal inequality for $\alpha$-mixing processes, due to Rio \cite{Rio:2017} Theorem 3.1 
\begin{align}
E\left(\max_{1\leq k \leq n} \left|\sum_{j=1}^{k}\xi_j \right| \right)^2 \leq 16\, \sum_{k=1}^n \int_0^1 \alpha^{-1}(u) \, Q_k^2(u)\, du, \label{rio}
\end{align}
where $Q_k$ denotes the quantile function of $\xi_k$, to the random variables $\xi_k=\frac{X_k}{\sqrt{k}}$. Note that
\begin{align*}
Q_k(u)=\frac{Q_X(u)}{\sqrt{k}}.
\end{align*}  
Furthermore, we use the inequality in Lemma \ref{ineq:SS}, i.e.\
\begin{align*}
\max_{1\leq k \leq n} \frac{\left|S_k\right|}{\sqrt{k}} \leq 2 \max_{1\leq k \leq n} \left| \sum_{j=1}^k \xi_j \right|.
\end{align*}
Then, for $K>0$, we obtain for all $n\geq n_K$

{ \allowdisplaybreaks
\begin{align*}
 & P\left(\sqrt{2\log\log n}\, Y_n -2\log\log n \geq -K\right)  \\
  \leq \, &  P\left(Y_n \geq \sqrt{\log\log n}\right) \\
  \leq \, &  P\left(2 \max_{ 1\leq k \leq (\log\log n)^\gamma }\left|\sum_{j=1}^k \xi_j \right| \geq \sqrt{\log\log n}    \right) \\
 \leq \, & \frac{4}{\log \log n} E\left(\max_{ 1\leq k \leq (\log\log n)^\gamma }\left|\sum_{j=1}^k \xi_j \right| \right)^2 \\
 \leq \, & \frac{4}{\log \log n} 16\, \sum_{k=1}^{(\log\log n)^\gamma} \int_0^1 \alpha^{-1}(u) \, Q_k^2(u)\, du \\
  = \, & \frac{4}{\log \log n} 16\, \sum_{k=1}^{(\log\log n)^\gamma} \int_0^1 \alpha^{-1}(u) \, \frac{Q_X^2(u)}{k}\, du \\
  \leq \, & \frac{4}{\log \log n} 16\, \sum_{k=1}^{(\log\log n)^\gamma} \int_0^1 \alpha^{-1}(u) \, \frac{Q_{|X|}^2(u)}{k}\, du \\
\leq \, &  C \frac{1}{\log\log n}\sum_{k=1}^{(\log \log n)^\gamma } \frac{1}{k} \\
\leq \, &  C \frac{\gamma \log\log\log n}{\log\log n} \longrightarrow 0.
\end{align*}
}
\end{proof}

\begin{lemma}
\label{lem:max-equiv}
Let $U_n$ and $V_n$ be real-valued random variables satisfying $U_n\cpr -\infty$. Then $V_n \claw G$ if and only if $\max(U_n,V_n)\claw G$. 
\end{lemma}
\begin{proof}
We have the following chain of inequalities
\begin{align*}
 P(V_n > x) \leq P(\max(U_n,V_n)>x) \leq P(U_n>x )+P(V_n>x),
\end{align*}
and hence $|P(\max(U_n,V_n) >x)-P(V_n>x)|\leq P(U_n>x) \rightarrow 0 $.
\end{proof}

\begin{proof}[Proof of Theorem~\ref{thm:DE}]

Let $(W)_{t\geq 0}$ be a Brownian motion with $\Var(W_1)=\sigma^2$. Define analogues of the random variables $Y_n$ and $Z_n$, replacing the partial sum process by Brownian motion
\begin{align*}
\widetilde{Y}_n&= \max_{1\leq k \leq (\log\log n)^\gamma} \frac{|W_k|}{\sqrt{k}} ,     \\
\widetilde{Z}_n&= \max_{(\log\log n)^\gamma\leq k \leq n} \frac{|W_k|}{\sqrt{k}}.  
\end{align*}
Note that $\widetilde{M}_n:=\max(\widetilde{Y}_n,\widetilde{Z}_n) = \max_{1\leq k \leq n} \frac{|W_k|}{\sqrt{k}} $. By the Darling-Erd\H{o}s theorem for Brownian motion, see Darling-Erd\H{o}s \cite{DE:1956}, we know that 
\begin{align*}
\frac{\sqrt{2\log\log n}}{\sigma} \; \widetilde{M}_n  -b_n \claw G.
\end{align*}
Applying Lemma~\ref{lem:y-n} to $\widetilde{Y}_n$, and Lemma~\ref{lem:max-equiv}, we obtain that 
\begin{align*}
 \frac{\sqrt{2\log\log n}}{\sigma} \; \widetilde{Z}_n  -b_n \claw G.
\end{align*}
In the final step, we employ an almost sure invariance principle which is stated as Theorem~1 in Merlev\`ede and Rio \cite{MR:2012} under a weaker strong mixing condition than we have. We can conclude that, under the assumptions of Theorem~\ref{thm:DE}, there exists a Brownian motion $(W_t)_{t\geq 0}$ with $\Var(W_1)=\sigma^2$, such that 
\begin{align*}
|S_k-W_k|=O\left(k^{\frac{1}{2}-\lambda}\right)~a.s.~\text{for}~ \lambda>0.
\end{align*}
Hence, we obtain
\begin{align*}
  \left| \max_{(\log\log n)^\gamma \leq k\leq n} \frac{|S_k|}{\sqrt{k}} -\max_{(\log\log n)^\gamma \leq k\leq n} 
\frac{|W_k|}{\sqrt{k}} \right| &\leq \max_{(\log\log n)^\gamma \leq k\leq n} \frac{|S_k-W_k|}{\sqrt{k}} \\
&\leq  C \max_{(\log\log n)^\gamma \leq k\leq n} \frac{k^{\frac{1}{2}-\lambda}}{\sqrt{k}} \\
&\leq  C (\log\log n)^{-\lambda \gamma}.
\end{align*}
Since $\lambda \gamma >\frac{1}{2} $, this implies that $\frac{\sqrt{\log\log n}}{\sigma} (\widetilde{Z}_n-Z_n)\cpr 0$. Hence, we obtain using Slutsky's lemma that
\begin{align*}
\frac{\sqrt{2\log\log n}}{\sigma} \; Z_n - b_n \claw G.
\end{align*}
Now, the statement of Theorem~\ref{thm:DE} follows from Lemma~\ref{lem:y-n} and Lemma~\ref{lem:max-equiv}.\\ 
\end{proof}

\appendix

\section{Auxiliary Results}
\begin{lemma}
\label{ineq:SS}
For all $a_1,\dots,a_k \in \mathbb{R}$, the following inequality holds 
\begin{align*}
\max_{1\leq k \leq n } \frac{1}{\sqrt{k}} \left|\sum_{j=1}^k a_j \sqrt{j}  \right| \leq 2 \max_{1\leq k \leq n} \left| \sum_{j=1}^k a_j \right|
\end{align*}
\end{lemma}
\begin{proof}
Define the partial sums $A_j:=\sum_{i=1}^j a_i$ and set $A_0=0$. Then $a_j=A_j-A_{j-1}$, and thus 
\begin{align*}
\left| \sum_{j=1}^k a_j \sqrt{j} \right| &= \left|\sum_{j=1}^k (A_j-A_{j-1})\sqrt{j} \right| \\
&= \left| \sum_{j=1}^{k-1} A_j \left(\sqrt{j}-\sqrt{j+1}\right)+A_k \sqrt{k} \right| \\
&\leq \max_{1\leq j \leq k} |A_j| \sum_{j=1}^{k-1} \left(\sqrt{j+1}-\sqrt{j}\right)+ |A_k|\sqrt{k} \\
&= \max_{1\leq j \leq k} |A_j| \sqrt{k} + |A_k| \sqrt{k} \\
&\leq 2 \sqrt{k} \max_{1\leq j \leq k } |A_j|.
\end{align*}
\end{proof}
\begin{remark}
This is a special case of an inequality stated as Lemma 1 in Shorack and Smythe \cite{SS:1976}.
\end{remark}

\begin{lemma}
\label{m.diss}
Assume that the kernel $g$ is degenerated and that it satisfies the variation condition. Let $m=\max\{i_{(2)}-i_{(1)},i_{(4)}-i_{(3)}\}$, where $\{i_1,i_2,i_3,i_4\}=\{i_{(1)},i_{(2)},i_{(3)},i_{(4)} \}$ and $i_{(1)}\leq i_{(2)} \leq i_{(3)} \leq i_{(4)}$. If $g$ is a kernel with uniform $(2+\delta)$-moments for a $\delta>0$, then there exists a constant $C$ such that 
\begin{align*}
\left|E\left(g(X_{i_1},X_{i_2})g(X_{i_3},X_{i_4})  \right)\right| \leq C \left( \int_{0}^{\alpha(m)} Q_{|X|}(u) du \right)^{\frac{\delta}{3+2\delta}}.
\end{align*}
\end{lemma}
\begin{proof}
Let $\varepsilon>0,~K>0$ and define
\begin{align*}
g_K(x,y)= \begin{cases} g(x,y) &~ \text{if}~|g(x,y)|\leq \sqrt{K} \\ \sqrt{K} &~\text{if}~ g(x,y)>\sqrt{K} \\ -\sqrt{K} & ~\text{if}~g(x,y)<-\sqrt{K}.  \end{cases}
\end{align*}
Let $g$ satisfy the variation condition with constant $L$. Then, $g_K$ satisfies also the variation condition with the same constant $L$. Assume, without loss of generality, that $m=i_2-i_1$. Moreover, we can assume that there exists a uniform on $[0,1]$ random variable that is independent of $(X_i)_{i\geq 1}$. With Theorem 1 of Peligrad \cite{P:2002}, choose a random  variable $X_{i_1}'$ independent of $X_{i_2},X_{i_3},X_{i_4}$ with the same distribution as $X_{i_1}$ and 
\begin{align*}
P\left(|X_{i_1}-X_{i_1}'|\geq \varepsilon \right) \leq \frac{4 \int_{0}^{\alpha(m)} Q_{|X|}(u) du}{\varepsilon}.
\end{align*}
As $g$ is a degenerate kernel, we have 
\begin{align*}
E\left(g(X_{i_1}',X_{i_2})g(X_{i_3},X_{i_4}) \right)=0.
\end{align*}
We get 
\begin{align}
& \left| E\left(g(X_{i_1},X_{i_2})g(X_{i_3},X_{i_4})  \right)  \right| \notag \\
=&  \left| E\left(g(X_{i_1},X_{i_2})g(X_{i_3},X_{i_4})  \right) - E\left(g(X_{i_1}',X_{i_2})g(X_{i_3},X_{i_4}) \right)  \right| \notag \\
=& E\left(\big| (g_K(X_{i_1},X_{i_2})-g_K(X_{i_1}',X_{i_2}))g_K(X_{i_3,X_{i_4}})\big| 1_{\{|X_{i_1}-X_{i_1}'|< \varepsilon \}}   \right) \label{delta1} \\
&+ E\left(\big| (g_K(X_{i_1},X_{i_2})-g_K(X_{i_1}',X_{i_2}))g_K(X_{i_3,X_{i_4}})\big| 1_{\{|X_{i_1}-X_{i_1}'|\geq \varepsilon \}}   \right)  \label{delta2} \\
&+ E\left(\big| g_K(X_{i_1},X_{i_2})g_K(X_{i_3},X_{i_4})-g(X_{i_1},X_{i_2})g(X_{i_3},X_{i_4})\big| \right)  \label{delta3} \\
&+ E\left(\big| g_K(X_{i_1}',X_{i_2})g_K(X_{i_3},X_{i_4})-g(X_{i_1}',X_{i_2})g(X_{i_3},X_{i_4})\big| \right). \label{delta4}
\end{align}
Due to the variation condition and $|g_K(X_{i_3},X_{i_4})|\leq \sqrt{K}$, which holds by definition, we have that \eqref{delta1} is smaller than $L\varepsilon\sqrt{K}$. The second summand \eqref{delta2} is bounded by
\begin{align*}
2KP(|X_{i_1}-X_{i_1}'|\geq \varepsilon)\leq \frac{8 \int_{0}^{\alpha(m)} Q_{|X|}(u) du}{\varepsilon} K.
\end{align*}
Let $M$ be the bound of the $(2+\delta)$-moments of $g$. For the third term of the sum \eqref{delta3} we get the following chain of inequalities.
\begin{align*}
& E\big(\big| g_K(X_{i_1},X_{i_2})g_K(X_{i_3},X_{i_4})-g(X_{i_1},X_{i_2})g(X_{i_3},X_{i_4})\big| \big) \\
\leq &  E \!\big(\! \big| g(X_{i_1},X_{i_2})\sqrt{K}\!-\!g(X_{i_1},X_{i_2})g(X_{i_3},X_{i_4})\big|\! 1_{\{|g(X_{i_1},X_{i_2})| \leq \sqrt{K}, g(X_{i_3},X_{i_4})> \sqrt{K} \}} \! \big) \\
& +\! E\!\big(\! \big|g(X_{i_1}\!,\!X_{i_2})\!(\!-\!\sqrt{K})\!-\!g(X_{i_1}\!,\!X_{i_2})g(X_{i_3}\!,\!X_{i_4})\big|\! 1_{\{\!|g(X_{i_1}\!,\!X_{i_2})| \leq\! \sqrt{K}, g(X_{i_3}\!,\!X_{i_4})<\!-\! \sqrt{K} \! \}} \! \big) \\
& + \! E\!\big( \!\big|\!\sqrt{K}g(X_{i_3},X_{i_4})\!-\!g(X_{i_1},X_{i_2})g(X_{i_3},X_{i_4})\big| \!1_{\{\!g(X_{i_1},X_{i_2}) > \sqrt{K}, |g(X_{i_3},X_{i_4})|\leq \sqrt{K}\! \}} \! \big) \\
& + \! E\!\big(\! \big|\!-\!\sqrt{K}g(X_{i_3}\!,\!X_{i_4})\!-\!g(X_{i_1}\!,\!X_{i_2})g(X_{i_3}\!,\!X_{i_4})\big|\! 1_{\{\!g(X_{i_1},X_{i_2})\!  < \!-\! \sqrt{K}, |g(X_{i_3},X_{i_4})|\leq \!\sqrt{K} \!\}} \! \big) \\
& + \! E\!\big(\! \big|\sqrt{K}\sqrt{K}-g(X_{i_1},X_{i_2})g(X_{i_3},X_{i_4})\big| 1_{\{g(X_{i_1},X_{i_2}) > \sqrt{K}, g(X_{i_3},X_{i_4})> \sqrt{K} \}} \! \big) \\
& + \! E\!\big( \!\big|\!(\!-\!\sqrt{K})(\!-\!\sqrt{K})\!-\!g(X_{i_1}\!,\!X_{i_2})g(X_{i_3}\!,\!X_{i_4})\big|\! 1_{\{\!g(X_{i_1},X_{i_2})< \!-\sqrt{K}, g(X_{i_3},X_{i_4})<-\! \sqrt{K} \! \}} \! \big) \\
& + \! E\!\big(\! \big|\sqrt{K}(-\sqrt{K})-g(X_{i_1},X_{i_2})g(X_{i_3},X_{i_4})\big|\! 1_{\{g(X_{i_1},X_{i_2})>\sqrt{K}, g(X_{i_3},X_{i_4})<- \sqrt{K} \}}  \!\big) \\
& + \! E\!\big( \!\big|(-\sqrt{K})\sqrt{K}-g(X_{i_1},X_{i_2})g(X_{i_3},X_{i_4})\big|\! 1_{\{g(X_{i_1},X_{i_2})< -\sqrt{K}, g(X_{i_3},X_{i_4})> \sqrt{K} \}} \! \big). 
\end{align*} 
Considering the first term of the sum, we have
\begin{align*}
& E\!\big( \!\big| g(X_{i_1},X_{i_2})\sqrt{K}\!-\!g(X_{i_1},X_{i_2})g(X_{i_3},X_{i_4})\big| 1_{\{|g(X_{i_1},X_{i_2})| \leq \sqrt{K}, g(X_{i_3},X_{i_4})> \sqrt{K} \}} \! \big) \\
= \, & E\big( \big|( \sqrt(K)- g(X_{i_3},X_{i_4}) ) g(X_{i_1},X_{i_2} \big| 1_{\{g(X_{i_3},X_{i_4})> \sqrt{K} \}}  \big) \\
\leq \, & \sqrt{K} E\big( \big| \sqrt(K)- g(X_{i_3},X_{i_4})\big| 1_{\{g(X_{i_3},X_{i_4})> \sqrt{K} \}}  \big) \\
\leq \, & \sqrt{K} E\big( (|g(X_{i_3},X_{i_4})|- \sqrt{K}) 1_{\{|g(X_{i_3},X_{i_4})|> \sqrt{K} \}}  \big) \\ 
\leq \, &  \sqrt{K} \frac{E|g(X_{i_3},X_{i_4})|^{2+\delta}}{(\sqrt{K})^{1+\delta}}  =  \frac{E|g(X_{i_3},X_{i_4})|^{2+\delta}}{K^{\delta/2}} \leq  \frac{M}{K^{\delta/2}}.
\end{align*}
The next three terms can be treated analogously. For the fifth summand we get the following inequality 
\begin{align*}
& E\big( \big|\sqrt{K}\sqrt{K}-g(X_{i_1},X_{i_2})g(X_{i_3},X_{i_4})\big| 1_{\{g(X_{i_1},X_{i_2}) > \sqrt{K}, g(X_{i_3},X_{i_4})> \sqrt{K} \}}  \big) \\
\leq  & E\big(\big( |g(X_{i_1},X_{i_2})| |g(X_{i_3},X_{i_4})|-K\big)1_{\{|g(X_{i_1},X_{i_2})| > \sqrt{K}, |g(X_{i_3},X_{i_4})|> \sqrt{K} \}}  \big) \\
\leq  & \!\big(\! E\!\big(\!|g(X_{i_1},X_{i_2})|^2 1_{\{|g(X_{i_1},X_{i_2})|>\sqrt{K}\}} \big) \! \big)^{1/2} \big(\! E\!\big(\!|g(X_{i_3},X_{i_4})|^2 1_{\{|g(X_{i_3},X_{i_4})|>\sqrt{K}\}} \!\big) \! \big)^{1/2} \\
\leq  & \left( \frac{E|g(X_{i_1},X_{i_2})|^{2+\delta}}{(\sqrt{K})^{\delta}}\right)^{1/2} \left( \frac{E|g(X_{i_3},X_{i_4})|^{2+\delta}}{(\sqrt{K})^{\delta}}\right)^{1/2} \leq  \frac{M}{K^{\delta/2}}.
\end{align*} 
Analogously, we get for the last three terms that they are also bounded by $\frac{M}{K^{\delta/2}}$. Altogether, we have
\begin{align*}
& E\big(\big| g_K(X_{i_1},X_{i_2})g_K(X_{i_3},X_{i_4})-g(X_{i_1},X_{i_2})g(X_{i_3},X_{i_4})\big| \big) \leq  8 \frac{M}{K^{\delta/2}}.
\end{align*}   
Similarly, we obtain for \eqref{delta4} 
\begin{align*}
E\big(\big| g_K(X_{i_1}',X_{i_2})g_K(X_{i_3},X_{i_4})-g(X_{i_1}',X_{i_2})g(X_{i_3},X_{i_4})\big| \big) \leq 8 \frac{M}{K^{\delta/2}}.
\end{align*}
Altogether, we get 
\begin{align*}
 \left| E\big(g(X_{i_1},X_{i_2})g(X_{i_3},X_{i_4})  \big)  \right| \leq  L \varepsilon \sqrt{K} + \frac{8 \int_{0}^{\alpha(m)} Q_{|X|}(u) du}{\varepsilon} K + 16\frac{M}{K^{\delta/2}}.
\end{align*}
Choosing $\varepsilon= \sqrt{ \int_{0}^{\alpha(m)} Q_{|X|}(u) du} K^{1/4}$ and $K=(\int_{0}^{\alpha(m)} Q_{|X|}(u) du)^{-\frac{1}{3/2+\delta}}$, we obtain
\begin{align*}
\left| E\big(g(X_{i_1},X_{i_2})g(X_{i_3},X_{i_4})  \big)  \right| \leq (L+8+16M) \left( \int_{0}^{\alpha(m)} Q_{|X|}(u) du \right)^{\frac{\delta}{3+2}}.
\end{align*}
\end{proof}

\begin{lemma}
\label{Lemma_DFGW}
In addition to the conditions of Lemma \ref{m.diss}, let \eqref{eq:mix-2} be satisfied. Then there exists a constant $C$ such that 
\begin{align*}
E\left( \sum_{i=1}^{k} \sum_{j=k+1}^{n} g(X_i,X_j)\right)^2 \leq C k(n-k). 
\end{align*}
\end{lemma}

\begin{proof}
This was proved for functionals of absolutely regular processes by Dehling et al.\ \cite{DFGW:2015} in Lemma~1. They make use of an upper bound for the expectations $|E(g(X_{i_1},X_{i_2})g(X_{i_3},X_{i_4}) )|$. Such a bound for an $\alpha$-mixing process is stated in Lemma \ref{m.diss}. The rest of the proof is analogous.
\end{proof}




\bibliographystyle{imsart-number} 
\bibliography{literature.bib}       

 

\end{document}